\documentclass[10pt]{article}
\usepackage{amsmath}
\usepackage{amssymb}
\usepackage{mathrsfs}
\usepackage{amsfonts}
\usepackage{mathrsfs,amscd,amssymb,amsthm,amsmath,bm,graphicx}
\usepackage{graphicx}
\usepackage{cite}
\usepackage{color}
\makeatletter

\renewcommand{\@seccntformat}[1]{{\csname the#1\endcsname}{\normalsize.}\hspace{.5em}}
\makeatother

\def \[{\begin{equation}}
\def \]{\end{equation}}

\newtheorem{thm}{Theorem}[section]

\newtheorem{lem}[thm]{Lemma}

\begin{document}
\setlength{\baselineskip}{13pt}
\begin{center}{\Large \bf The Laplacian and normalized Laplacian spectra of M\"{o}bius polyomino networks and their applications}

\vspace{4mm}

{\large Zhi-Yu Shi $^{1}$, Jia-Bao Liu $^{1,2,*}$, Sakander Hayat $^{3}$}\vspace{2mm}

{\small $^{1}$School of Mathematics and Physics, Anhui Jianzhu University, Hefei 230601, China\\
$^{2}$School of Mathematics, Southeast University, Nanjing 210096, China\\
$^{3}$Faculty of Engineering Sciences, GIK Institute of Engineering Sciences and Technology, Topi 23460,
Pakistan}
\vspace{2mm}
\end{center}

\footnotetext{E-mail address: shizhiyuah@163.com, liujiabaoad@163.com, sakander1566@gmail.com}

\footnotetext{* Corresponding author.}

{\noindent{\bf Abstract.}\ \ Spectral theory has widely used in complex networks and solved some practical problems.
In this paper, we investigated the Laplacian and normalized Laplacian spectra of M\"{o}bius polyomino networks by using spectral theory. Let $M_{n}$ denote M\"{o}bius polyomino networks ($n\geq3$). As applications of the obtained results, the Kirchhoff index, multiplicative degree-Kirchhoff index, Kemeny's constant and spanning trees of $M_{n}$ are obtained. Moreover, it is surprising to find that the multiplicative degree-Kirchhoff index of $M_{n}$ is nine times as much as the Kirchhoff index.

 \vspace{2mm}

\section{Introduction}\label{sct1}
\ \ \ \ \ In 1964, Heilbronner \cite{A} proposed M\"{o}bius aromatic based on Huckel's molecular orbital theory. Compared with Huckel's system, the M\"{o}bius system is stable because of its closed shell structure. In recent years, compounds with M\"{o}bius aromatic have been synthesized.
In particular, Ma et al. \cite{B} studied the normalized Laplacian spectrum for the hexagonal M\"{o}bius graphs. Then, Geng et al. \cite{C} obtained the Laplacian spectrum of M\"{o}bius phenylenes Chain. In 2019, Lei et al. \cite{1} studied the normalized Laplacian spectrum of M\"{o}bius phenylene chain.
For other graphs, see \cite{2,3,4,5,6,7,8}. Motivated by these, we investigate the Laplacian spectrum and normalized Laplacian spectrum of M\"{o}bius polyomino networks. The M\"{o}bius polyomino networks are one kind molecular graphs embedded into the M\"{o}bius strip, which each side is a polyomino (see Figure 1(1)). For example, the M\"{o}bius polyomino network of length 6 is depicted in Figure 1(2).

\begin{figure}[htbp]
\centering\includegraphics[width=10cm,height=8cm]{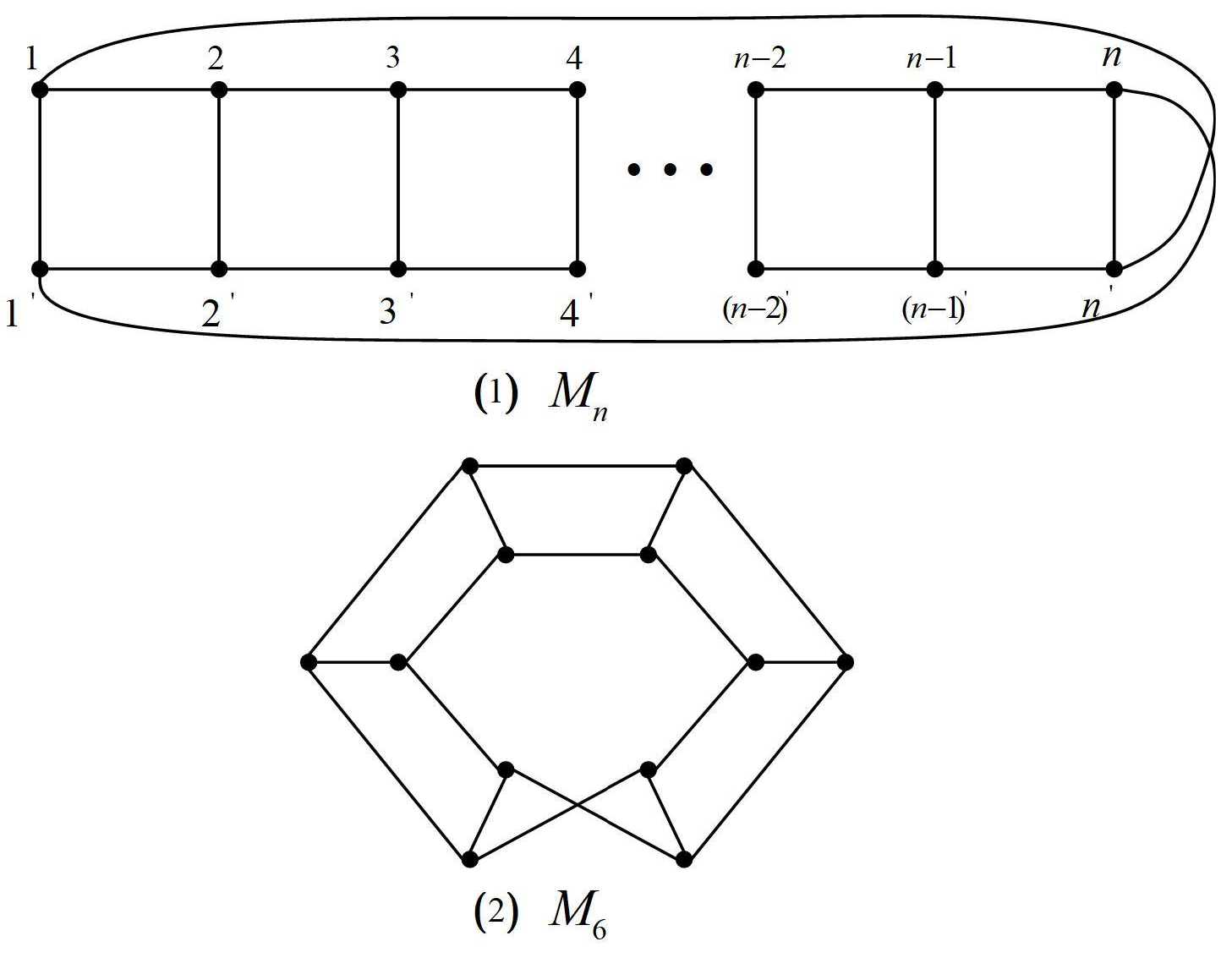}
\caption{M\"{o}bius polyomino networks.}
\end{figure}

Since $n=1$, the graph does not exist. When $n=2$, the graph corresponds to the special graph in reference \cite{D}. Therefore, we investigate the results of M\"{o}bius polyomino networks which $n\geq 3$.

Let $G=(E_{G},V_{G})$ be a graph with edge set $E_{G}$ and vertex set $V_{G}$, where its size are $|E_G|=m$ and $|V_G|=n$. Let $D(G)=diag\{d_{1},d_{2},\cdots,d_{n}\}$ represent a degree matrix, and $A(G)$ be the adjacency matrix, where $d_{i}$ is the degree of $v_{i}$. The Laplacian matrix of $G$ is given by $L(G)=D(G)-A(G)$, which $(i,j)$-entry are equal to -1 when $v_i$ and $v_j$ are adjacent, $d_i$ when $i=j$, 0 otherwise. The normalized Laplacian matrix of $G$ is given by $\mathcal{L}(G)=D(G)^{-\frac{1}{2}}LD(G)^{-\frac{1}{2}}$, which $(i,j)$-entry are equal to $-\frac{1}{\sqrt{d_id_j}}$ when $v_i$ and $v_j$ are adjacent, 1 when $i=j$, 0 otherwise. For more notation, one can be referred to \cite{E}.

The traditional concept of distance is the length of the shortest path between vertices $i$ and $j$, represented by $d_{ij}$. The Wiener index \cite{a} was proposed as
$W(G)=\sum_{i<j}d_{ij}$.
Wiener index has been extensively studied in chemistry. In the past two decades, the research on the Wiener index are shown in references \cite{b,c,d,e}.
In 1994, the Gutman index \cite{f} was proposed as
$Gut(G)=\sum_{i<j}d_{i}d_{j}d_{ij}$.

Klein and Randi\'c \cite{F} was the first to put forward the concept of resistance distance, and the resistance distance between vertices $i$ and $j$ is denoted by $r_{ij}$.  One famous resistance distance-based parameter called the Kirchhoff index \cite{F},
namely $Kf(G)=\sum_{i<j}r_{ij}$. In 2007, Chen et al. \cite{G} defined the multiplicative degree-Kirchhoff index as $Kf^{*}(G)=\sum_{i<j}d_{i}d_{j}r_{ij}$. For convenience, Gutman and Zhu et al. \cite{H,I} introduced the Kirchhoff index as
\begin{eqnarray}
Kf(G)=n\sum_{k=2}^{n}\frac{1}{\mu_{k}},
\end{eqnarray}
where $0=\mu_{1}<\mu_{2}\leq\cdots\leq\mu_{n}$ are the Laplacian eigenvalues of $G$.

Chen et al. \cite{G} proposed the multiplicative degree-Kirchhoff index as
\begin{eqnarray}
Kf^{*}(G)=2m\sum_{k=2}^{n}\frac{1}{\lambda_{k}},
\end{eqnarray}
where $0=\lambda_{1}<\lambda_{2}\leq\cdots\leq\lambda{n}$ are the normalized Laplacian eigenvalues of $G$.

In Section \ref{sct2}, we mainly introduce some notations and theorems. Then, the Laplacian spectrum of $M_{n}$ is investigated in Section \ref{sct3}. In Section \ref{sct4}, we obtained the normalized Laplacian spectrum of $M_{n}$ in the same way as in Section \ref{sct3}. The conclusion and discussion are summarized in Section \ref{sct5}.

\section{Preliminary}\label{sct2}
\ \ \ \ \ Given an $n\times n$ matrix $B$, submatrix of $B$ is represented by $B[i_{1},\cdots,i_{k}]$, where $B[i_{1},\cdots,i_{k}]$ is formed by removing the $i_{1}$-th,$\cdots$,$i_{k}$-th rows and columns of $B$. Let $P_{B}(x)=det(xI-B)$ represent characteristic polynomial of $B$.

Label M\"{o}bius polyomino networks as shown in the Figure 1(1). Evidently, $|V(M_{n})|=2n,|E(M_{n})|=3n$ and $V_{1}=\{1,2,\cdots,n\},V_{2}=\{1',2',\cdots,n'\}$ is an automorphism of $M_{n}$.

Then $L(M_n)$ and $\mathcal{L}(M_n)$ can be expressed by

\begin{equation*}
L(M_{n})=\left(
\begin{array}{cc}
L_{V_{1}V_{1}}& L_{V_{1}V_{2}}\\
L_{V_{2}V_{1}}& L_{V_{2}V_{2}}\\
\end{array}
\right),~
\mathcal{L}(M_{n})=\left(
\begin{array}{cc}
\mathcal{L}_{V_{1}V_{1}}& \mathcal{L}_{V_{1}V_{2}}\\
\mathcal{L}_{V_{2}V_{1}}& \mathcal{L}_{V_{2}V_{2}}\\
\end{array}
\right),
\end{equation*}
where
\begin{equation*}
L_{V_{1}V_{1}}=L_{V_{2}V_{2}},~ L_{V_{1}V_{2}}=L_{V_{2}V_{1}},~ \mathcal{L}_{V_{1}V_{1}}=\mathcal{L}_{V_{2}V_{2}},~ \mathcal{L}_{V_{1}V_{2}}=\mathcal{L}_{V_{2}V_{1}}.
\end{equation*}

Let
\begin{equation*}
T=\left(
  \begin{array}{cc}
  \frac{1}{\sqrt{2}}I_{n}& \frac{1}{\sqrt{2}}I_{n}\\
  \frac{1}{\sqrt{2}}I_{n}& -\frac{1}{\sqrt{2}}I_{n}
  \end{array}
\right),
\end{equation*}
then
\begin{equation*}
TL(M_{n})T=\left(
  \begin{array}{cc}
  L_{A}& 0\\
    0& L_{S}
  \end{array}
\right),~
T\mathcal{L}(M_{n})T=\left(
  \begin{array}{cc}
  \mathcal{L}_{A}& 0\\
    0& \mathcal{L}_{S}
  \end{array}
\right),
\end{equation*}
where
\begin{eqnarray*}
L_{A}=L_{V_{1}V_{1}}+L_{V_{1}V_{2}},
L_{S}=L_{V_{1}V_{1}}-L_{V_{1}V_{2}},
\mathcal{L}_{A}=\mathcal{L}_{V_{1}V_{1}}+\mathcal{L}_{V_{1}V_{2}},
\mathcal{L}_{S}=\mathcal{L}_{V_{1}V_{1}}-\mathcal{L}_{V_{1}V_{2}}.
\end{eqnarray*}

In what follows, the theorems that we present will be used throughout the Section 3 and Section 4.

\begin{thm}\textup{\cite{J}}
If $L_{A},L_{S},\mathcal{L}_{A},\mathcal{L}_{S}$ are defined as above, the following formula can be obtained. Then
\begin{eqnarray*}
P_{L(M_{n})}(x)=P_{L_{A}}(x)P_{L_{S}}(x), P_{\mathcal{L}(M_{n})}(x)=P_{\mathcal{L}_{A}}(x)P_{\mathcal{L}_{S}}(x).
\end{eqnarray*}
\end{thm}

\begin{thm}\textup{\cite{K,L}}
For a cycle with n vertices $C_n$:

(1) The Kirchhoff index of $C_n$ is
\begin{eqnarray*}
Kf(C_n)=\frac{n^{3}-n}{12}.
\end{eqnarray*}

(2) The Laplacian eigenvalues of $C_n$ is
\begin{eqnarray*}
\alpha_{i}=2-2\cos\frac{2\pi i}{n},i\in[1,n].
\end{eqnarray*}
\end{thm}

\begin{thm}\textup{\cite{M}}
For a connected graph $G$ with vertices $n$. The number of spanning trees of $G$ is
\begin{eqnarray*}
\tau(G)=\frac{1}{n}\prod_{i=2}^{n}\mu_{i},
\end{eqnarray*}
where $\mu_{i}$ is the Laplacian eigenvalue of $G$.
\end{thm}

\begin{thm}\textup{\cite{N}}
For a connected graph $G$ with vertices $n$. Kemeny's constant of $G$ is
\begin{eqnarray*}
Kc(G)=\sum_{i=2}^{n}\frac{1}{\lambda_i},
\end{eqnarray*}
where $\lambda_{i}$ is the normalized Laplacian eigenvalue of $G$.
\end{thm}

Evidently, from Theorem 2.4 and (1.4), the relation between the Kemeny's constant and multiplicative degree-Kirchhoff index is $Kf^{*}(G)=2mKc(G)$.

\begin{thm}\textup{\cite{M}}
If graph $G$ with $|V_G|=n$ and $|E_G|=m$. The number of spanning trees of $G$ is
\begin{eqnarray*}
2m\tau(G)=\prod_{i=1}^{n}d_i\cdot\prod_{i=2}^{n}\lambda_i,
\end{eqnarray*}
where $\lambda_{i}$ is the normalized Laplacian eigenvalue of $G$.
\end{thm}

For convenience, all $p$ represents $2+\sqrt{3}$ and all $q$ represents $2-\sqrt{3}$ throughout the paper.

\section{Laplacian spectrum of $M_n$}\label{sct3}
\ \ \ \ \ In this section, we mainly obtain the Laplacian spectrum of $M_n$.

According to Laplacian matrix of $M_n$, we can get $L_{V_1 V_1}$ and $L_{V_1 V_2}$:

\begin{eqnarray*}
 L_{V_1 V_1}&=&
\left(
  \begin{array}{cccccc}
    3 & -1 & & & &\\
    -1 & 3 & -1 & & &\\
    & -1 & 3 & -1 & &\\
    & & & \ddots & &\\
    & & & -1 & 3 & -1\\
    & & & & -1 & 3\\
  \end{array}
\right)_{n\times n},
\end{eqnarray*}

\begin{eqnarray*}
 L_{V_1 V_2}&=&
\left(
  \begin{array}{cccccc}
    -1 & & & & & -1\\
    & -1 & & & &\\
    & & -1 & & &\\
    & & & \ddots & &\\
    & & & & -1 &\\
    -1 & & & & & -1\\
  \end{array}
\right)_{n\times n}.
\end{eqnarray*}

Based on Theorem 2.1, Laplacian spectrum consists of the eigenvalues of $L_{A}$ and $L_{S}$ of $M_{n}$ can be obtained.

\begin{eqnarray*}
 L_{A}&=&
\left(
  \begin{array}{cccccc}
    2 & -1 & & & & -1\\
    -1 & 2 & -1 & & &\\
    & -1 & 2 & -1 & &\\
    & & & \ddots & &\\
    & & & -1 & 2 & -1\\
    -1 & & & & -1 & 2\\
  \end{array}
\right)_{n\times n},
\end{eqnarray*}

\begin{eqnarray*}
 L_{S}&=&
\left(
  \begin{array}{cccccc}
    4 & -1 & & & & 1\\
    -1 & 4 & -1 & & &\\
    & -1 & 4 & -1 & &\\
    & & & \ddots & &\\
    & & & -1 & 4 & -1\\
    1 & & & & -1 & 4\\
  \end{array}
\right)_{n\times n}.
\end{eqnarray*}

Assume that
$0=\alpha_{1}<\alpha_{2}\leq\alpha_{3}\leq\cdots\leq\alpha_{n}$ are the
roots of $P_{L_{A}}(x)=0$, and
$0<\beta_{1}\leq\beta_{2}\leq\beta_{3}\leq\cdots\leq\beta_{n}$ are the
roots of $P_{L_{S}}(x)=0$. Noticing that $L_A$ is the Laplacian matrix of cycle $C_{n}$, the following lemma can be obtained from (1.1) and Theorem 2.2(1).

\begin{lem}
For M\"{o}bius polyomino networks $M_n$ ($n\geq3$),
\begin{eqnarray*}
Kf(M_n)=\frac{n^{3}-n}{6}+2n\sum_{i=1}^{n}\frac{1}{\beta_{i}},
\end{eqnarray*}
where $\beta_{i}$ is the eigenvalue of $L_S$.
\end{lem}

Next, we first determine $\sum_{i=1}^{n}\frac{1}{\beta_{i}}$.

\begin{lem}
Let $0<\beta_{1}\leq\beta_{2}\leq\beta_{3}\leq\cdots\leq\beta_{n}$ be the eigenvalues of $L_S$. Then
\begin{eqnarray*}
\sum_{i=1}^{n}\frac{1}{\beta_{i}}=\frac{n}{2\sqrt{3}}\cdot\frac{p^{n}-q^{n}}{p^{n}+q^{n}+2}.
\end{eqnarray*}
\end{lem}

\noindent{\bf Proof.} Let
\begin{eqnarray*}
P_{L_{S}}(x)=det(xI-L_{S})=x^{n}+a_{1}x^{n-1}+\cdots+a_{n-1}x+a_{n}.
\end{eqnarray*}

According to Vieta$^{'}$s theorem of $P_{L_{S}}(x)$, one obtain
\begin{eqnarray}
\sum_{i=1}^{n}\frac{1}{\beta_{i}}=
\frac{(-1)^{n-1}a_{n-1}}{det(L_S)}.
\end{eqnarray}

Since $(-1)^{n-1}a_{n-1}$ and $det(L_S)$, we focus on $i$-th order principal submatrix $F_i$, which consists of the first $i$ rows and columns of the following matrix $L_{S}^{'}$, $i\in[1,n]$. Let $f_i=det(F_i)$.

\begin{eqnarray*}
 L_{S}^{'}&=&
\left(
  \begin{array}{cccccc}
    4 & -1 & & & &\\
    -1 & 4 & -1 & & &\\
    & -1 & 4 & -1 & &\\
    & & & \ddots & &\\
    & & & -1 & 4 & -1\\
    & & & & -1 & 4\\
  \end{array}
\right)_{n\times n}.
\end{eqnarray*}

It is easy to get that $f_1=4,f_2=15,f_3=56$ and for $3\leq{i}\leq{n}$,
$f_i=4f_{i-1}-f_{i-2}$. For convenience, we let $f_0=1$. Furthermore, one can verity that
\begin{eqnarray*}
f_i=\frac{p^{i+1}-q^{i+1}}{2\sqrt{3}}.
\end{eqnarray*}

\noindent{\bf Fact 1.} $(-1)^{n-1}a_{n-1}=\frac{n}{2\sqrt{3}}(p^{n}-q^{n}).$

\noindent{\bf Proof of Fact 1.} Obviously, we obtain $(-1)^{n-1}a_{n-1}$ is the sum of all the principal minors of order $n-1$ of $L_S$. According to the property of $L_S$, we know that
\begin{eqnarray*}
det(L_S[i])=
\begin{cases}
f_{n-1}, & i=1~or~n;\\
f_{i-1}f_{n-i}-f_{i-2}f_{n-i-1}, & i\in[2,n-1].
\end{cases}
\end{eqnarray*}

Thus, we can obtain
\begin{eqnarray*}
(-1)^{n-1}a_{n-1}&=&\sum_{i=1}^{n}det(L_S[i])\\
&=&det(L_S[1])+det(L_S[n])+\sum_{i=2}^{n-1}det(L_S[i])\\
&=&\frac{n}{2\sqrt{3}}(p^{n}-q^{n}).
\end{eqnarray*}

This result as desired. \hfill\rule{1ex}{1ex}

\noindent{\bf Fact 2.} $det(L_S)=p^{n}+q^{n}+2.$

\noindent{\bf Proof of Fact 2.} By expanding the last row of $L_{S}$, we can arrive at
\begin{eqnarray*}
det(L_S)&=&f_n-f_{n-2}+2\\
&=&p^{n}+q^{n}+2.
\end{eqnarray*}

Thus Fact 2 holds. \hfill\rule{1ex}{1ex}

Substituting Facts 1 and 2 into (3.3) yields lemma 3.2. \hfill\rule{1ex}{1ex}

In view of lemmas 3.1 and 3.2, the following theorem can be obtained.

\begin{thm}
For M\"{o}bius polyomino networks $M_n$ ($n\geq3$),
\begin{eqnarray*}
Kf(M_n)=\frac{n^{3}-n}{6}+\frac{n^{2}}{\sqrt{3}}\cdot\frac{p^{n}-q^{n}}{p^{n}+q^{n}+2}.
\end{eqnarray*}
\end{thm}

\noindent{\bf Remark.} We found that Theorem 3.3 can be obtained by different methods in \cite{P} and \cite{Q}.

In view of the Laplacian spectrum of $M_n$, the following theorem can be obtained.

\begin{thm}
For M\"{o}bius polyomino networks $M_n$ ($n\geq3$),
\begin{eqnarray*}
\tau(M_n)=\frac{n}{2}(p^{n}+q^{n}+2).
\end{eqnarray*}
\end{thm}

\noindent{\bf Proof.} Based on Theorem 2.2(2), we know that the eigenvalues of $L_A$ are $\alpha_{i}=2-2\cos\frac{2\pi{i}}{n}$ ($1\leq{i}\leq{n}$) and the product of $\alpha_{2},\alpha_{3},\cdots,\alpha_{n}$ are $\prod_{i=2}^{n}\alpha_{i}=n^{2}$. By Fact 2, we get
\begin{eqnarray*}
\prod_{i=1}^{n}\beta_{i}=det(L_S)=p^{n}+q^{n}+2.
\end{eqnarray*}

Together with Fact 2 and Theorem 2.3, Theorem 3.4 follows immediately. \hfill\rule{1ex}{1ex}

\section{Normalized Laplacian spectrum of $M_n$}\label{sct4}
\ \ \ \ \ In this section, we mainly obtain the normalized Laplacian spectrum of $M_n$.
Moreover, we find that the multiplicative degree-Kirchhoff index of $M_n$ is nine times of Kirchhoff index.

According to normalized Laplacian matrix of $M_n$, we can get $\mathcal{L}_{V_1 V_1}$ and $\mathcal{L}_{V_1 V_2}$:

\begin{eqnarray*}
\mathcal{L}_{V_1 V_1}&=&
\left(
  \begin{array}{cccccc}
    1 & -\frac{1}{3} & & & &\\
    -\frac{1}{3} & 1 & -\frac{1}{3} & & &\\
    & -\frac{1}{3} & 1 & -\frac{1}{3} & &\\
    & & & \ddots & &\\
    & & & -\frac{1}{3} & 1 & -\frac{1}{3}\\
    & & & & -\frac{1}{3} & 1\\
  \end{array}
\right)_{n\times n},
\end{eqnarray*}

\begin{eqnarray*}
\mathcal{L}_{V_1 V_2}&=&
\left(
  \begin{array}{cccccc}
    -\frac{1}{3} & & & & & -\frac{1}{3}\\
    & -\frac{1}{3} & & & &\\
    & & -\frac{1}{3} & & &\\
    & & & \ddots & &\\
    & & & & -\frac{1}{3} &\\
    -\frac{1}{3} & & & & & -\frac{1}{3}\\
  \end{array}
\right)_{n\times n}.
\end{eqnarray*}

Based on Theorem 2.1, we can get normalized Laplacian spectrum which is obtained by the eigenvalues of $\mathcal{L}_{A}$ and $\mathcal{L}_{S}$ of $M_{n}$.

\begin{eqnarray*}
\mathcal{L}_{A}&=&
\left(
  \begin{array}{cccccc}
    \frac{2}{3} & -\frac{1}{3} & & & & -\frac{1}{3}\\
    -\frac{1}{3} & \frac{2}{3} & -\frac{1}{3} & & &\\
    & -\frac{1}{3} & \frac{2}{3} & -\frac{1}{3} & &\\
    & & & \ddots & &\\
    & & & -\frac{1}{3} & \frac{2}{3} & -\frac{1}{3}\\
    -\frac{1}{3} & & & & -\frac{1}{3} & \frac{2}{3}\\
  \end{array}
\right)_{n\times n},
\end{eqnarray*}

\begin{eqnarray*}
\mathcal{L}_{S}&=&
\left(
  \begin{array}{cccccc}
    \frac{4}{3} & -\frac{1}{3} & & & & \frac{1}{3}\\
    -\frac{1}{3} & \frac{4}{3} & -\frac{1}{3} & & &\\
    & -\frac{1}{3} & \frac{4}{3} & -\frac{1}{3} & &\\
    & & & \ddots & &\\
    & & & -\frac{1}{3} & \frac{4}{3} & -\frac{1}{3}\\
    \frac{1}{3} & & & & -\frac{1}{3} & \frac{4}{3}\\
  \end{array}
\right)_{n\times n}.
\end{eqnarray*}

Assume that
$0=\gamma_{1}<\gamma_{2}\leq\gamma_{3}\leq\cdots\leq\gamma_{n}$ are the
roots of $P_{\mathcal{L}_{A}}(x)=0$, and
$0<\delta_{1}\leq\delta_{2}\leq\delta_{3}\leq\cdots\leq\delta_{n}$ are the
roots of $P_{\mathcal{L}_{S}}(x)=0$. Next, we derive the formulas of $\sum_{i=2}^{n}\frac{1}{\gamma_{i}}$ and $\sum_{i=1}^{n}\frac{1}{\delta_{i}}$.

\begin{lem}
Let $0=\gamma_{1}<\gamma_{2}\leq\gamma_{3}\leq\cdots\leq\gamma_{n}$ be the eigenvalues of $\mathcal{L}_A$. Then
\begin{eqnarray*}
\sum_{i=2}^{n}\frac{1}{\gamma_{i}}=\frac{n^{2}-1}{4}.
\end{eqnarray*}
\end{lem}

\noindent{\bf Proof.} Suppose that
\begin{eqnarray*}
P_{\mathcal{L}_{A}}(x)=det(xI-\mathcal{L}_{A})=x^{n}+b_{1}x^{n-1}+\cdots+b_{n-2}x^{2}+b_{n-1}x.
\end{eqnarray*}

Applying Vieta$^{'}$s theorem, one can get
\begin{eqnarray}
\sum_{i=2}^{n}\frac{1}{\gamma_{i}}=\frac{(-1)^{n-2}b_{n-2}}{(-1)^{n-1}b_{n-1}}.
\end{eqnarray}

Before calculating $(-1)^{n-1}b_{n-1}$ and $(-1)^{n-2}b_{n-2}$, we focus on $i$-th order principal submatrix $G_i$, which consists of the first $i$ rows and columns of the following matrix $\mathcal{L}_{A}^{'}$, $i\in[1,n]$. Let $g_i=det(G_i)$.

\begin{eqnarray*}
\mathcal{L}_{A}^{'}&=&
\left(
  \begin{array}{cccccc}
    \frac{2}{3} & -\frac{1}{3} & & & &\\
    -\frac{1}{3} & \frac{2}{3} & -\frac{1}{3} & & &\\
    & -\frac{1}{3} & \frac{2}{3} & -\frac{1}{3} & &\\
    & & & \ddots & &\\
    & & & -\frac{1}{3} & \frac{2}{3} & -\frac{1}{3}\\
    & & & & -\frac{1}{3} & \frac{2}{3}\\
  \end{array}
\right)_{n\times n}.
\end{eqnarray*}

It is easy to get that $g_1=\frac{2}{3},g_2=\frac{1}{3},g_3=\frac{4}{27}$ and for $3\leq{i}\leq{n}$,
$g_i=\frac{2}{3}g_{i-1}-\frac{1}{9}g_{i-2}$. For convenience, we let $g_0=1$. Furthermore, we can get
\begin{eqnarray*}
g_i=\frac{i+1}{3^{i}}.
\end{eqnarray*}

\noindent{\bf Fact 3.} $(-1)^{n-1}b_{n-1}=\frac{n^{2}}{3^{n-1}}.$

\noindent{\bf Proof of Fact 3.} Evidently, we obtain $(-1)^{n-1}b_{n-1}$ is the sum of all the principal minors of order $n-1$ of $\mathcal{L}_A$. According to the property of $\mathcal{L}_A$, we know that
\begin{eqnarray*}
det(\mathcal{L}_A[i])=
\begin{cases}
g_{n-1}, & i=1~or~n;\\
g_{i-1}g_{n-i}-\frac{1}{9}g_{i-2}g_{n-i-1}, & i\in[2,n-1].
\end{cases}
\end{eqnarray*}

Therefore, one can get
\begin{eqnarray*}
(-1)^{n-1}b_{n-1}&=&\sum_{i=1}^{n}det(\mathcal{L}_{A}[i])\\
&=&det(\mathcal{L}_{A}[1])+det(\mathcal{L}_{A}[n])+\sum_{i=2}^{n-1}det(\mathcal{L}_{A}[i])\\
&=&\frac{n^{2}}{3^{n-1}},
\end{eqnarray*}
which is the desired result. \hfill\rule{1ex}{1ex}

\noindent{\bf Fact 4.} $(-1)^{n-2}b_{n-2}=\frac{n^{2}(n^{2}-1)}{4\cdot3^{n-1}}.$

\noindent{\bf Proof of Fact 4.} Obviously, we obtain $(-1)^{n-2}b_{n-2}$ is the sum of all the principal minors of order $n-2$ of $\mathcal{L}_A$. In a similar way, we know that
\begin{eqnarray*}
det(\mathcal{L}_A[i,j])=
\begin{cases}
g_{n-2}, & i=1~and~j=n;\\
g_{i-2}g_{n-i}, & j=1~or~n,~and~i\in[2,n-1];\\
g_{i-1}g_{j-i-1}g_{n-j}-\frac{1}{9}g_{i-2}g_{j-i-1}g_{n-j-1}, & i<j~and~i,j\in[2,n-1].
\end{cases}
\end{eqnarray*}

Thus, one can obtain
\begin{eqnarray*}
(-1)^{n-2}b_{n-2}&=&\sum_{1\leq{i}<j\leq{n}}det(\mathcal{L}_A[i,j])\\
&=&det(\mathcal{L}_A[1,n])+\sum_{i=2}^{n-1}det(\mathcal{L}_A[1,i])+\sum_{i=2}^{n-1}det(\mathcal{L}_A[i,n])+\sum_{2\leq{i}<j\leq{n-1}}det(\mathcal{L}_A[i,j])\\
&=&\frac{n^{2}(n^{2}-1)}{4\cdot3^{n-1}}.
\end{eqnarray*}

This completes the proof. \hfill\rule{1ex}{1ex}

Substituting Facts 3 and 4 into (4.4) yields lemma 4.1. \hfill\rule{1ex}{1ex}

\begin{lem}
Let $0<\delta_{1}\leq\delta_{2}\leq\delta_{3}\leq\cdots\leq\delta_{n}$ be the eigenvalues of $\mathcal{L}_S$. Then
\begin{eqnarray*}
\sum_{i=1}^{n}\frac{1}{\delta_{i}}=\frac{\sqrt{3}n}{2}\cdot\frac{p^{n}-q^{n}}{p^{n}+q^{n}+2}.
\end{eqnarray*}
\end{lem}

\noindent{\bf Proof.} Let
\begin{eqnarray*}
P_{\mathcal{L}_{S}}(x)=det(xI-\mathcal{L}_{S})=x^{n}+k_{1}x^{n-1}+\cdots+k_{n-1}x+k_n.
\end{eqnarray*}

Applying Vieta$^{'}$s theorem, one can get
\begin{eqnarray}
\sum_{i=1}^{n}\frac{1}{\delta_{i}}=\frac{(-1)^{n-1}k_{n-1}}{det(\mathcal{L}_S)}.
\end{eqnarray}

In order to obtain $(-1)^{n-1}k_{n-1}$ and $det(\mathcal{L}_S)$, we focus on $i$-th order principal submatrix $H_i$, which consists of the first $i$ rows and columns of the following matrix $\mathcal{L}_{S}^{'}$, $i\in[1,n]$. Let $h_i=det(H_i)$.

\begin{eqnarray*}
\mathcal{L}_{S}^{'}&=&
\left(
  \begin{array}{cccccc}
    \frac{4}{3} & -\frac{1}{3} & & & &\\
    -\frac{1}{3} & \frac{4}{3} & -\frac{1}{3} & & &\\
    & -\frac{1}{3} & \frac{4}{3} & -\frac{1}{3} & &\\
    & & & \ddots & &\\
    & & & -\frac{1}{3} & \frac{4}{3} & -\frac{1}{3}\\
    & & & & -\frac{1}{3} & \frac{4}{3}\\
  \end{array}
\right)_{n\times n}.
\end{eqnarray*}

It is easy to get that $h_1=\frac{4}{3},h_2=\frac{5}{3},h_3=\frac{56}{27}$ and for $3\leq{i}\leq{n}$,
$h_i=\frac{4}{3}h_{i-1}-\frac{1}{9}h_{i-2}$. For convenience, we let $h_0=1$. Furthermore, one can verity that
\begin{eqnarray*}
h_i=\frac{p^{i+1}-q^{i+1}}{2\sqrt{3}\cdot3^{i}}.
\end{eqnarray*}

\noindent{\bf Fact 5.} $(-1)^{n-1}k_{n-1}=\frac{n}{2\sqrt{3}}\cdot\frac{p^{n}-q^{n}}{3^{n-1}}.$

\noindent{\bf Proof of Fact 5.} Similarly, we obtain $(-1)^{n-1}k_{n-1}$ is the sum of all the principal minors of order $n-1$ of $\mathcal{L}_S$. Thus, we know that
\begin{eqnarray*}
det(\mathcal{L}_S[i])=
\begin{cases}
h_{n-1}, & i=1~or~n;\\
h_{i-1}h_{n-i}-\frac{1}{9}h_{i-2}h_{n-i-1}, & i\in[2,n-1].
\end{cases}
\end{eqnarray*}

Therefore, one can obtain
\begin{eqnarray*}
(-1)^{n-1}k_{n-1}&=&\sum_{i=1}^{n}det(\mathcal{L}_S[i])\\
&=&det(\mathcal{L}_{S}[1])+det(\mathcal{L}_{S}[n])+\sum_{i=2}^{n-1}det(\mathcal{L}_S[i])\\
&=&\frac{n}{2\sqrt{3}}\cdot\frac{p^{n}-q^{n}}{3^{n-1}}.
\end{eqnarray*}

This completes the proof. \hfill\rule{1ex}{1ex}

\noindent{\bf Fact 6.} $det(\mathcal{L}_S)=\frac{p^{n}+q^{n}+2}{3^{n}}.$

\noindent{\bf Proof of Fact 6.} By expanding the last row of $\mathcal{L}_{S}$, we can arrive at
\begin{eqnarray*}
det(\mathcal{L}_S)&=&h_n-\frac{1}{9}h_{n-2}+\frac{2}{3^{n}}\\
&=&\frac{p^{n}+q^{n}+2}{3^{n}}.
\end{eqnarray*}

Thus Fact 6 holds. \hfill\rule{1ex}{1ex}

Substituting Facts 5 and 6 into (4.5) yields lemma 4.2. \hfill\rule{1ex}{1ex}

Based on Theorem 2.4, lemmas 4.1 and 4.2, the following theorem can be obtained.

\begin{thm}
For M\"{o}bius polyomino networks $M_n$ ($n\geq3$),
\begin{eqnarray*}
Kc(M_n)=\frac{n^{2}-1}{4}+\frac{\sqrt{3}n}{2}\cdot\frac{p^{n}-q^{n}}{p^{n}+q^{n}+2}.
\end{eqnarray*}
\end{thm}

In view of (1.2), lemmas 4.1 and 4.2, one get the multiplicative degree-Kirchhoff index.

\begin{thm}
For M\"{o}bius polyomino networks $M_n$ ($n\geq3$),
\begin{eqnarray*}
Kf^{*}(M_n)=\frac{3}{2}(n^{3}-n)+\frac{3\sqrt{3}n^{2}(p^{n}-q^{n})}{p^{n}+q^{n}+2}.
\end{eqnarray*}
\end{thm}

According to Theorem 3.3 and Theorem 4.4, we find that the multiplicative degree-Kirchhoff index of $M_n$ is nine times of Kirchhoff index. Since the degree of each point in the M\"{o}bius polyomino networks is three, $Kf^{*}(M_n)=d_id_jKf(M_n)=9Kf(M_n)$ can be obtained from the definition of Kirchhoff index and multiplicative degree-Kirchhoff index index. Therefore, it is not difficult to obtain and prove the correctness of results.

In view of the normalized Laplacian spectrum of $M_n$, the following theorem can be obtained.

\begin{thm}
For M\"{o}bius polyomino networks $M_n$ ($n\geq3$),
\begin{eqnarray*}
\tau(M_n)=\frac{n}{2}(p^{n}+q^{n}+2).
\end{eqnarray*}
\end{thm}

\noindent{\bf Proof.} Based on Theorem 2.5, one can get $\prod_{i=1}^{2n}d_i\prod_{i=2}^{n}\frac{1}{\gamma_{i}}\prod_{i=1}^{n}\frac{1}{\delta_{i}}=2\cdot{3n}\cdot\tau(M_n)$,
where

$$\prod_{i=1}^{2n}d_i=3^{2n},$$

$$\prod_{i=2}^{n}\frac{1}{\gamma_{i}}=(-1)^{n-1}b_{n-1}=\frac{n^{2}}{3^{n-1}},$$

$$\prod_{i=1}^{n}\frac{1}{\delta_{i}}=det(\mathcal{L}_S)=\frac{p^{n}+q^{n}+2}{3^{n}}.$$

Hence,

\begin{eqnarray*}
\tau(M_n)=\frac{n}{2}(p^{n}+q^{n}+2).
\end{eqnarray*}

The result as desired. \hfill\rule{1ex}{1ex}

\section{Conclusion and discussion}\label{sct5}
\ \ \ \ \ In this paper, based on the Laplacian matrix and normalized Laplacian matrix of M\"{o}bius polyomino networks, the Kirchhoff index, multiplicative degree-Kirchhoff index, Kemeny's constant and spanning trees of M\"{o}bius polyomino networks are determined through the decomposition theorem and Vieta$^{'}$s theorem.

Spectral theory has important applications in many fields. Wu \cite{7} and Ding \cite{O} obtained the mean first-passage time of Koch networks and 3-prism graph according to the Laplacian spectrum, respectively. Thus, we can explore the mean first-passage time of M\"{o}bius polyomino networks.

\section*{Funding}
\ \ \ \ \ This work was supported in part by National Natural Science Foundation of China Grant 11601006. 

\end{document}